\documentclass[12pt]{article}
\usepackage[russian]{babel}
\usepackage{amsfonts}
\usepackage{amssymb,amsmath,amsthm,eufrak,euscript}
\usepackage{graphicx,fancybox}

\newtheorem{lemma}{Lemma}
\newtheorem{theorem}{Theorem}
\newtheorem{example}{Example}
\begin{document}

\begin{center}
{\Large Integrable magnetic geodesic flows on 2-torus: new example via
quasi-linear system of PDEs.}
\end{center}

\medskip

\begin{center}
{\large S.V. Agapov\footnote[1]{Sobolev Institute of Mathematics,
Novosibirsk, Russia. Supported by RSF (grant 14-11-00441)}, M.
Bialy\footnote[2]{School of Mathematical Sciences, Raymond and
Beverly Sackler Faculty of Exact Sciences, Tel Aviv University,
Israel. Supported by ISF grant 162/15.}, A.E. Mironov$^1$}
\end{center}

\begin{quote}
\noindent{\sc Abstract. } The only one example has been known of
magnetic geodesic flow on the 2-torus which has a polynomial in
momenta integral independent of the Hamiltonian. In this example the
integral is linear in momenta and corresponds to a one parametric
group preserving the Lagrangian function of the magnetic flow. In
this paper the problem of integrability on one energy level is
considered. This problem can be reduced to a remarkable Semi-hamiltonian system of
quasi-linear PDEs and to the question of existence of smooth periodic
solutions for this system. Our main result states that the pair of
Liouville metric with zero magnetic field on the 2-torus can be
analytically deformed to a Riemannian metric with small magnetic
field so that the magnetic geodesic flow on an energy level is
integrable by means of a quadratic in momenta integral. Thus our
construction gives a new example of smooth periodic solution to the
Semi-hamiltonian (Rich) quasi-linear system of PDEs.

\medskip

\end{quote}

\section{\large Introduction and main results}
Consider  Hamiltonian system with two degrees of freedom
$$
\dot{x}^j = \{x^j,H\}, \qquad \dot{p}_j = \{p_j,H\}, \qquad j=1,2
\eqno
$$
with a Hamiltonian function $H$. Let us remind, that the Hamiltonian
system is called integrable if there exists a function (the first
integral) $F(x_1,x_2,p_1,p_2),$ such that
$$
\dot{F} = \{F,H\} = \sum_{i=1}^2 \left ( \frac{\partial F}{\partial
x^i} \frac{\partial H}{\partial p_i} - \frac{\partial F}{\partial
p_i} \frac{\partial H} {\partial x^i} \right ) = 0,
$$
and $F$ is functionally independent with $H$ almost everywhere.

Geodesic flow of a Riemannian metric $ds^2 = g_{ij}dx^idx^j$ on the
2-torus can be viewed as a Hamiltonian system with Hamiltonian
function $H = \frac{1}{2} g^{ij} p_ip_j.$ On the 2-torus there are
two known types of metrics with integrable geodesic flow. Namely, in
the first case the metric in some coordinates has the form
$$ds^2 = \Lambda (x^1) ((dx^1)^2 + (dx^2)^2).$$
In this case there exists an additional integral which is linear in
momenta
$$F = p_2.$$
In the second case the metric is Liouville, that is in some
coordinates it has the form
$$ds^2 = (\Lambda_1 (x^1) + \Lambda_2 (x^2)) ((dx^1)^2 +
(dx^2)^2).$$ Then there exists an integral quadratic in momenta
$$F = \frac{\Lambda_2 (x^2) p_1^2-\Lambda_1 (x^1) p_2^2 }
{\Lambda_1 (x^1) + \Lambda_2 (x^2)}. \eqno(1)$$ The question of
existence of metrics on the 2-torus with polynomial first integral
of higher degree is old and turns out to be very complicated
(see~\cite{5},~\cite{6}, ~\cite{BMq},\-~\cite{14},~\cite{12},~\cite{15}).

The question goes back to G. Birkhoff~\cite{Birkhoff}, who completed the study
integrals of degree 1 and 2, and in more modern time it was
conjectured in~\cite{BFK} that the two cases
described above exhaust all examples of integrable geodesic flows on
the 2-torus. This conjecture is not known even for polynomial
integrals of degrees 3,4. In our recent papers (see~\cite{4}--~\cite{BMq})
we invented new approach to the problem and
proved in particular that it is equivalent to the problem of
existence of smooth periodic solutions of a system of quasi-linear
PDEs of the form (3) (see below) which turns out to be semi-Hamiltonian (or
Rich). Moreover in some cases using~\cite{18}, the
behavior of smooth solutions and their blow up can be understood for
these quasi-linear systems~\cite{3},~\cite{BMq},~\cite{BMe}.

Hamiltonian system of magnetic geodesic flow is described by the
same Hamiltonian function $H = \frac{1}{2} g^{ij} p_ip_j,$ but with
the symplectic structure twisted by the 2-form of the magnetic field
$\Omega=\Omega(x_1,x_2)dx_1\wedge dx_2.$ Thus in the magnetic case
the Poisson bracket takes the form:
$$
\{F,H\}_{mg} = $$ $$=\sum_{i=1}^2 \left ( \frac{\partial F}{\partial
x^i} \frac{\partial H}{\partial p_i} - \frac{\partial F}{\partial
p_i} \frac{\partial H}{\partial x^i} \right ) + \Omega (x^1,x^2)
\left ( \frac{\partial F}{\partial p_1} \frac{\partial H}{\partial
p_2} - \frac{\partial F}{\partial p_2} \frac{\partial H}{\partial
p_1} \right ). \eqno(2)
$$
Many questions related to integrability of magnetic flows have been
considered in the literature (~\cite{2},~\cite{7},~\cite{8},~\cite{9},~\cite{E},~\cite{20},~\cite{21}). The following, is the only one known example
 of polynomially integrable magnetic geodesic
flows on 2-torus with non-zero magnetic field:
\begin{example}
Let the Riemannian metric be of the form
$$ds^2 = \Lambda (y) (dx^2 + dy^2),$$
and the magnetic form $$\Omega=-u'(y)dx\wedge dy .$$ Then the first
integral is linear in momenta and takes the form
$$F_1 = p_1 + u (y).$$
\end{example}
It is plausible that no other examples of integrable magnetic flows
exist on the 2-torus. We shall prove this guess in Section 4 for
quadratic in momenta integrals. The underlying idea here is very
natural: the dynamics of magnetic flow is changing if the energy
levels of $H$ varies. So the requirement for $F$ to be the first
integral on all levels is a heavy requirement. However, the proof of
this fact in the general case is not known and might be complicated.

Therefore, we turn now to the much more interesting problem of
existence of Polynomial integral on a particular level
$\{H=\frac{1}{2}\}.$ It turns out, that the problem of existence of
integrable magnetic flows on one energy level can be reduced to a
fundamental question on existence of smooth
 periodic solutions of a quasi-linear PDEs of the following form:
$$
A(U) U_x + B(U) U_y = 0, \eqno(3)
$$
where the matrices $A(U), B(U)$ depend explicitly on the components
of the vector $U = (u_1, \ldots, u_n)^T$ containing the metric
factor $\Lambda$ and some coefficients of the polynomial integral
$F$.

The system of the form (3) appears for integrable geodesic flows and
magnetic geodesic flows. It was derived in~\cite{4} --~\cite{6} and proved
that (3) has infinitely many conservation laws and in the hyperbolic
region it can be written in the form of Riemann invariants. Such
systems are called Semi-hamiltonian systems of hydrodynamic type (or
Rich). They were introduced in~\cite{22},~\cite{23} (see also in~\cite{18}).
We will discuss the properties of (3) in Section 2 in more details.

Let us notice that if the magnetic field is equal to zero
identically, then for $n=2$, the system (3) is linearly degenerate
and this property of the system enables the existence of smooth
periodic solutions which give rise to the first integral of the form
(1). When the magnetic field is present, it turns out that system
(3) is not linearly degenerate anymore and moreover the sign of the
derivatives $\frac{\partial\lambda_i}{\partial r_i}$ changes along
characteristics, where $r_i$ is the Riemann invariant, $\lambda_i$ is the corresponding eigenvalue, so the question of existence of smooth periodic
solutions becomes very complicated. The main challenge here is to
prove or disprove the existence of periodic solutions of these
systems. The main result of the present paper is that for the case
of quadratic integral, i.e. $n=2$, and the system is $4\times 4$, we
managed to overcome these difficulties to get:

\begin{theorem}
There exist real analytic Riemannian metrics on the 2-torus which
are arbitrary close to the Liouville metrics (and different from
them) and a non-zero analytic magnetic fields such that magnetic
geodesic flows on the energy level $\{H=\frac{1}{2}\}$ have
polynomial in momenta first integral of degree two.
\end{theorem}

Let us stress again that our theorem gives a very natural $4\times
4$ Rich quasi-linear system which is not linearly degenerate having
analytic periodic solutions. We refer to~\cite{10} for
other examples of a very different nature.

The idea of the proof of the main theorem is as follows. Using the
specific form of the matrices of the system (3) for $n=2$, we
introduce another quasi-linear system of evolution equations
(formula (11) below)  which form a \emph{symmetry} of our system (3).
The flow of this system acts on the space of solutions. Moreover, we
choose the Liouville metric together with the zero magnetic field as
the initial data of this flow. Then this flow for any small $t>0$
generates the nontrivial Riemannian metrics and a non-zero magnetic
field which satisfy the requirements of Theorem 1. It is remarkable
fact, which is not hard to verify, that the set of ''symmetric``
solutions with respect to this symmetry consists only of flat
metrics with zero magnetic fields so does not produce any
interesting example.

\section{\large Semi-hamiltonian systems of hydrodynamic type}

Systems like (3) arise when geodesic flows are studied.
Such systems have many beautiful properties.

It is shown in~\cite{5} that in the case of the integrable geodesic flow one may introduce global
semi-geodesic coordinates $(t,x)$ on the 2-torus such that
$$
ds^2 = g^2(t,x)dt^2+dx^2, \qquad H = \frac{1}{2}
\left( \frac{p_1^2}{g^2}+p_2^2 \right).
$$

The condition $\dot{F} = \{F,H\} = 0,$ where the first integral $F$ has the form
$$
F = \frac{a_0}{g^n} p_1^n + \frac{a_1}{g^{n-1}} p_1^{n-1}p_2 + \ldots + \frac{a_{n-2}}{g^2} p_1^2 p_2^{n-2} + \frac{a_{n-1}}{g}p_1p_2^{n-1} + a_np_2^n, \ a_k = a_k(t,x)
$$
is equivalent to the system of quasi-linear PDEs of the form (3) with $U = (a_0, \ldots, a_{n-2},
a_{n-1})^T, \ a_{n-1}=g, \ a_n=1.$ In this case $A$ is the identity matrix and
\[ B = \left( \begin{array}{ccccccccc}
0 & 0 & \ldots & 0 & 0 & a_1\\
a_{n-1} & 0 & \ldots & 0 & 0 & 2a_2-na_0\\
0 & a_{n-1} & \ldots & 0 & 0 & 3a_3-(n-1)a_1\\
\ldots & \ldots & \ldots & \ldots & \ldots & \ldots\\
0 & 0 & \ldots & a_{n-1} & 0 & (n-1)a_{n-1}-3a_{n-3}\\
0 & 0 & \ldots & 0 & a_{n-1} & na_n-2a_{n-2}
\end{array} \right). \]
This system has many interesting properties. It can be
written in the form of conservation laws, i.e. there exists a change
of variables $ U^T \rightarrow (Q_1(U), \ldots, Q_n(U)) $ such that
for some $R_1(U), \ldots, R_n(U)$ the following relati\-ons hold:
$$
(Q_j(U))_t + (R_j(U))_x = 0, \qquad j =1, \ldots, n.
$$
Moreover, in the hyperbolic region, where all the eigenvalues
$\lambda_1, \ldots, \lambda_n$ of matrix $B$ are real and pairwise
distinct the system possesses $n$ Riemann's invariants
$
r_1(U), \ldots, r_n(U)
$
such that the system can be written in the following form:
$$
(r_j)_t + \lambda_j(r)(r_j)_x = 0, \qquad j =1, \ldots, n.
$$
Such systems are called semi-hamiltonian (see~\cite{5}). Semi-hamiltonian systems were introduced and studied by S.P. Tsarev in~\cite{22},~\cite{23} (see also~\cite{18},~\cite{19}).

The same question but in conformal coordinates $ds^2 =
\Lambda (x,y) (dx^2+dy^2)$ was studied in~\cite{6}. Given a polynomial
integral of the form
$$
F = a_0 p_1^n + a_1 p_1^{n-1} p_2 + \ldots +a_n p_2^n, \qquad a_k =
a_k(x,y),
$$
 it then follows from Kolokoltsov's theorem (see~\cite{12}) that the following
relations hold true
$$
a_n = c_1+a_{n-2}-a_{n-4}+ \ldots, \qquad a_{n-1} = c_2+a_{n-3}-a_{n-5}+ \ldots,
$$
where $c_1, c_2$ are some constants. The rest of unknown
coefficients and the conformal factor satisfy the quasi-linear
system of PDEs of the form (3), where $U = (a_0, \ldots, a_{n-2},
\Lambda)^T$. In this case this system also turns out to be
semi-hamiltonian (in the regions where at least one of the matrixes
$A$ and $B$ is non-degenerate).

For any semi-hamiltonian system the following relations on the eigenvalues hold:
$$
\partial_{r_j} \frac{\partial_{r_i} \lambda_k}{\lambda_i-\lambda_k} =
\partial_{r_i} \frac{\partial_{r_j} \lambda_k}{\lambda_j-\lambda_k},
\qquad i \neq j \neq k \neq i.
$$
It means that there is a diagonal metric
$$
ds^2 = H_1^2 (r) dr_1^2 + \ldots + H_n^2 (r) dr_n^2 \eqno(4)
$$
with Christoffel symbols satisfying the following relations:
$$
\Gamma_{ki}^k = \frac{\partial_{r_i} \lambda_k}{\lambda_i-\lambda_k}, \qquad i \neq k.
$$

It is proved in~\cite{6} that in the case of conformal coordinates the rotation coefficients $\beta_{kl}$ of the metric (4) associated with the system (3) are symmetric:
$$
\beta_{kl} = \beta_{lk}, \qquad \beta_{kl} = \frac{\partial_{r_k} H_l}{H_k}, \qquad k \neq l.
$$
This means that there is a function $a(r)$ such that $\partial_{r_k}
a(r) = H_k^2(r).$ Here $H_i$ are Lame coefficients of the metric
(4), $H_i^2 = g_{ii}.$ Such metrics are called the metrics of Egorov
type. Following~\cite{17}, we shall call the corresponding
semi-hamiltonian systems to be of Egorov type.

As it is shown in~\cite{17}, if the system does
not split
$(\partial_{r_i} \lambda_k \not\equiv 0, \ i \neq k)$, then it is of
Egorov type iff it possesses two special conservation laws:
$$
P_x+Q_y=0, \qquad P_y+R_x=0.
$$
In the case of conformal coordinates such conservation laws were found explicitly for the system (3) in~\cite{6}.

Similar results hold for the magnetic geodesic flow on the 2-torus.
Such flows (or, equiva\-lently, systems with gyroscopic forces) were
studied, for example, in~\cite{7},~\cite{13},~\cite{20},~\cite{21}. In conformal coordinates the
existence of an additional polynomial in momenta first integral of
an arbitrary degree $N$ on the fixed energy level $\{H =
\frac{1}{2}\}$ leads to the PDEs of the form (3) which is also
proved to be semi-hamiltonian (see~\cite{4}). Moreover, in~\cite{4} it's
proved that in the case of $N=2,3$ this is the Egorov type system.
Recently this result was generalized to the case of an arbitrary
degree $N$ (see~\cite{1}).

\section{The proof of the theorem}

Consider Hamiltonian system
$$
\dot{x}^j = \{x^j,H\}_{mg}, \qquad \dot{p}_j =
 \{p_j,H\}_{mg}, \qquad j=1,2 \eqno(5)
$$
on a 2-torus in the presence of magnetic field with Hamiltonian $H =
\frac{1}{2} g^{ij} p_ip_j$ and the Poisson bracket of the form (2).
Choose the conformal coordinates $(x,y)$ in which
$$ds^2 = \Lambda (x,y) (dx^2+dy^2), \qquad H = \frac{p_1^2+p_2^2}{2 \Lambda}.$$
Then one can parametrize the momenta on the energy level
$\{H=\frac{1}{2}\}$ in the following way:
$$
p_1 = \sqrt{\Lambda} \cos \varphi, \qquad p_2 = \sqrt{\Lambda} \sin \varphi.
$$
The equations (5) take the form
$$
\dot{x} = \frac{\cos \varphi}{\sqrt{\Lambda}},
\qquad \dot{y} = \frac{\sin \varphi}{\sqrt{\Lambda}},
\qquad \dot{\varphi} = \frac{\Lambda_y}{2 \Lambda \sqrt{\Lambda}} \cos \varphi -
\frac{\Lambda_x}{2 \Lambda \sqrt{\Lambda}} \sin \varphi - \frac{\Omega}{\Lambda}.
$$

Following~\cite{4}, we shall search for the first integral $F$ of the
second degree in the following way:
$$
F(x,y,\varphi) = \sum_{k=-2}^{k=2} a_k(x,y) e^{i k \varphi},\quad
a_k= u_k+ i v_k,\quad a_{-k} = \bar{a}_k. \eqno(6)
$$
The condition $\dot{F}=0$ is equivalent to the following
equation:
$$
F_x \cos \varphi + F_y \sin \varphi +
F_{\varphi} \left( \frac{\Lambda_y}{2 \Lambda} \cos \varphi -
\frac{\Lambda_x}{2 \Lambda} \sin \varphi - \frac{\Omega}{\sqrt{\Lambda}} \right)
= 0. \eqno(7)
$$

Let's substitute (6) into (7) and equate the coefficients at
$e^{i k \varphi}$ to zero. We obtain
$$
\frac{\Lambda_y}{2 \Lambda} \frac{i (k-1) a_{k-1} + i (k+1) a_{k+1}}{2} -
\frac{\Lambda_x}{2 \Lambda} \frac{i (k-1) a_{k-1} - i (k+1) a_{k+1}}{2 i} +
$$
$$
+\frac{(a_{k-1})_x+(a_{k+1})_x}{2} + \frac{(a_{k-1})_y-(a_{k+1})_y}{2 i} -
\frac{i k \Omega a_k}{\sqrt{\Lambda}} = 0, \eqno(8)
$$
where $k = 0, \ldots, 3,$ and $ \ a_k$ are assumed to be zero when
$k > 2.$ For $k=3$, we get from (8):
$$
 (a_2)_x-\frac{\Lambda_x}{\Lambda}a_2+\frac{1}{i}\left((a_2)_y-
 \frac{\Lambda_y}{\Lambda}a_2\right)=0.
$$
Multiplying this identity by $\Lambda^{-1}$ we obtain:
$$
 (a_2\Lambda^{-1})_x-i(a_2\Lambda^{-1})_y=0.
$$
Thus $a_2\Lambda^{-1}$ is anti-holomorphic function on 2-torus and
hence must be constant.  Consequently  $a_2=(\alpha+i\beta)\Lambda$
for some non-zero constant $\alpha+i\beta$. Applying appropriate
rotation in the plane $(x,y)$ and dividing $F$ by appropriate
constant we can assume that $\alpha=1,\beta=0$. Thus we get
$a_2=\Lambda$. Moreover, multiplying by $\Lambda^{-\frac{1}{2}}$ the
imaginary part of (8) for $k=2$ we compute:
$$
\Omega=\frac{1}{4}\left( (v_1\Lambda^{-\frac{1}{2}})_x-
(u_1\Lambda^{-\frac{1}{2}})_y  \right).
$$
Let us introduce new functions:
$$
f = \frac{u_1}{\sqrt{\Lambda}}, \qquad g = \frac{v_1}{\sqrt{\Lambda}}.
$$
Then, we get the following system of equations on the unknowns
$\Lambda, u_0,f,g$:
$$
 f_x+g_y=0,
$$
$$
 (f\Lambda)_x-(g\Lambda)_y=0,
$$
$$
 (u_0)_x+2\Lambda_x-\frac{1}{2}g(f_y-g_x)=0,
$$
$$
 -(u_0)_y+2\Lambda_y+\frac{1}{2}f(f_y-g_x)=0,
$$
which is equivalent to the matrix form:
$$
A(U)U_x+B(U)U_y=0, \eqno(9)
$$
where
\[
A = \begin{pmatrix} 0 & 0 & 1 & 0\\
f & 0 & \Lambda & 0\\
2 & 1 & 0 & \frac{g}{2}\\
0 & 0 & 0 & -\frac{f}{2} \end{pmatrix}, \qquad
B = \begin{pmatrix} 0 & 0 & 0 & 1\\
-g & 0 & 0 & -\Lambda\\
0 & 0 & -\frac{g}{2} & 0\\
2 & -1 & \frac{f}{2} & 0 \end{pmatrix},
\]
$U = (\Lambda, u_0,f,g)^T.$ The magnetic field takes the form
$$\Omega = \frac{1}{4}(g_x-f_y), \eqno(10)$$
hence the magnetic form is exact. One can check that
\[
U_0(x,y) = \begin{pmatrix} \Lambda_1(x)+\Lambda_2(y)\\
2 \Lambda_2(y)-2\Lambda_1(x)\\
0\\
0\end{pmatrix}
\]
is the solution of the system (9), where $\Lambda_1(x)$ and
$\Lambda_2(y)$ are periodic positive functions: $\Lambda_1(x+1) =
\Lambda_1(x), \ \Lambda_2(y+1) = \Lambda_2(y).$ This solution
corresponds to the case of geodesic flow of the Liouville metric
with zero magnetic field having the first integral (1). In the
sequel we assume $\Lambda_1$ and $\Lambda_2$ to be real analytic
functions.

We introduce the following evolution equations which are explicitly
written in the form:
$$
U_t=A_1(U)U_x+B_1(U)U_y, \eqno(11)
$$
where
\[
A_1 = \begin{pmatrix} g & 0 & 0 & \Lambda\\
-2g & g & 0 & -2 \Lambda\\
0 & 0 & 0 & 0\\
0 & -2 & 0 & 0 \end{pmatrix}, \qquad
B_1 = \begin{pmatrix} f & 0 & \Lambda & 0\\
2f & f & 2 \Lambda & 0\\
0 & 2 & 0 & 0\\
0 & 0 & 0 & 0 \end{pmatrix}.
\]
It turns out that this system of evolution equations define the
symmetry of the system (9) so that the flow of (11) transforms
solutions to solutions as we shall prove below.

Next we apply the following consequence of Cauchy--Kowalevskaya
theorem:

\begin{lemma}
The Cauchy problem for the system (11) with the initial data
$$U(x,y,t)\mid_{t=0} = U_0(x,y)\eqno(12)$$
has a unique analytic periodic $(U(x+1,y,t) = U(x,y+1,t) =
U(x,y,t))$ solution for $t$ small enough.
\end{lemma}
\textit{ Proof. } Since the entries of the matrices $A(U), B(U)$ are
linear functions in the components of $U$, it follows from Cauchy--Kowalevskaya theorem (see~\cite{11}) that there exists an analytic
solution $U(x,y,t)$ of the system (11) defined on a domain $D$
containing the plane $\{t=0\}$. Moreover, since the initial data of
(11) is chosen to be 1-periodic in $x,y$ it follows from the
uniqueness part of the theorem that the domain $D$ can be chosen
invariant under the 1-shifts $x\rightarrow x+1$ and $y\rightarrow
y+1$ and the solution $U(x,y,t)$ is 1-periodic in $x,y$. This proves
Lemma 1.

Let us prove that $U(x,y,t)$ constructed in Lemma 1 is a solution of our system
(9) for all small $t$. We denote by $\tilde{V}(x,y,t)$ the following
real analytic vector function
$$
\tilde{V} = A(U) U_x + B(U) U_y.
$$
By our construction $\tilde{V}(x,y,0) = 0.$ We have to prove that
$\tilde{V} \equiv 0.$ Denote
\[
\tilde{V} = \begin{pmatrix} V_1\\
V_2\\
V_3\\
V_4\end{pmatrix}, \qquad
V = \begin{pmatrix} V_2\\
V_3\\
V_4\end{pmatrix}.
\]
By direct calculations using (11) one can check that $\tilde{V}$
satisfies the following system of equations:
\[
\begin{pmatrix} V_1\\
V_2\\
V_3\\
V_4\end{pmatrix}_t= A_2\begin{pmatrix} V_1\\
V_2\\
V_3\\
V_4\end{pmatrix}_x+B_2\begin{pmatrix} V_1\\
V_2\\
V_3\\
V_4\end{pmatrix}_y+C_2\begin{pmatrix} V_1\\
V_2\\
V_3\\
V_4\end{pmatrix}+D_2\begin{pmatrix} V_1^2\\
V_1V_2\\
V_1V_3\\
V_2V_3\end{pmatrix}. \eqno(13)
\]
Here
\[
A_2 = \begin{pmatrix} 0 & 0 & 0 & 0\\
-g \Lambda & g & 0 & -2 \Lambda\\
0 & 0 & 0 & 0\\
2 \Lambda & -2 & f & g \end{pmatrix}, \qquad
B_2 = \begin{pmatrix} 0 & 0 & 0 & 0\\
f \Lambda & f & 2 \Lambda & 0\\
2 \Lambda & 2 & f & g\\
0 & 0 & 0 & 0 \end{pmatrix},
\]
\[
C_2 = \begin{pmatrix} 0 & 0 & 0 & 0\\
c_{1} & c_{2} & c_{3} & -2 \Lambda_x\\
c_{4} & 0 & f_y & -f_x\\
c_{5} & 0 & -g_y & g_x \end{pmatrix}, \qquad
D_2 = \begin{pmatrix} 0 & 0 & 0 & 0\\
-\frac{f \Lambda}{g} & -\frac{f}{g} & -\frac{4 \Lambda}{g} & -\frac{4}{g}\\
0 & 0 & 0 & 0\\
0 & 0 & 0 & 0 \end{pmatrix},
\]
where
$$
c_{1}=\frac{2 \Lambda f_x f + (f^2 - g^2 + 4 \Lambda) \Lambda_x +
2 \Lambda u_{0_x}}{g}, \ c_{2}=\frac{2 (g g_x + 2 \Lambda_x + u_{0_x})}{g},
$$
$$
c_{3}=\frac{4 \Lambda f_x +
2 f \Lambda_x}{g}, \ c_{4}= 4 \Lambda_y + \frac{1}{2} f (f_y-g_x),
\ c_{5}= 4 \Lambda_x - \frac{1}{2} g (f_y-g_x).
$$
So we have
$$
V_{1_t} = 0,
$$
hence
$$
V_1(x,y,t) = 0.
$$
Taking this into account, we can simplify the system (13). Let us
write it down explicitly:
$$
V_{2_t} = \frac{2 (g g_x + 2 \Lambda_x + u_{0_x})}{g} V_2 +
\frac{4 \Lambda f_x + 2 f \Lambda_x}{g} V_3  -
\frac{4}{g} V_2V_3 + g V_{2_x} -2 \Lambda V_{4_x} +
$$
$$
+f V_{2_y} + 2 \Lambda V_{3_y} - 2 \Lambda_x V_4, \eqno(14)
$$
$$
V_{3_t} = 2 V_{2_y} + f V_{3_y} + g V_{4_y} +f_y V_3 - f_x V_4,
$$
$$
V_{4_t} = -2 V_{2_x} + f V_{3_x} +g V_{4_x} - g_y V_3 +g_x V_4.
$$
Let us remark here that $g\equiv 0$ for ${t=0}$ and therefore Cauchy-Kowalevskaya theorem does not apply directly. But notice that for $t>0,$ $g$
becomes different from $0$, since from the evolution equations we
have $g_t=4\Lambda_{1x}\neq 0,$ for all but finitely many $x$. Next we prove:

\begin{lemma}
The Cauchy problem for the system (14) with the initial data
$$V(x,y,t)\mid_{t=0} = 0$$ has the unique analytic solution $V \equiv 0$.
\end{lemma}

\textit{Proof.} Let us assume that we have analytical solution $V$
of this Cauchy problem. Expand $V$ as well as $U$ in the power
series in $t$:

\[
V(x,y,t) = \sum_{k=0}^{+\infty} W_{k}(x,y) t^k, \ W_{k}(x,y) =\begin{pmatrix}
w_{2k}\\
w_{3k}\\
w_{4k}\end{pmatrix}, \ W_{0}(x,y) =\begin{pmatrix}
0\\
0\\
0\end{pmatrix}, \eqno(15)
\]

$$
\Lambda (x,y,t) = \sum_{k=0}^{+\infty} \lambda_{k}(x,y) t^k, \qquad \lambda_{0}(x,y) = \Lambda_1(x) + \Lambda_2(y), \eqno(16)
$$
$$
u_0(x,y,t) = \sum_{k=0}^{+\infty} r_k(x,y) t^k, \qquad r_0(x,y) = 2 \Lambda_2(y) - 2 \Lambda_1(x),
$$
$$
f(x,y,t) = \sum_{k=0}^{+\infty} s_{k}(x,y) t^k, \qquad s_{0}(x,y) = 0,
$$
$$
g(x,y,t) = \sum_{k=0}^{+\infty} q_{k}(x,y) t^k, \qquad q_{0}(x,y) = 0.
$$
Let us substitute these expansions into (14). We obtain the following identities:
$$
t q_1 w_{21} + O(t^2) = 0,
$$
$$
w_{31} + t (2 w_{32}-2 w_{{21}_y}) + O(t^2) =0,
$$
$$
w_{41} + t (2 w_{42}+2 w_{{21}_x}) + O(t^2) =0,
$$
it then follows that $w_{31} = 0, \ w_{41} = 0.$ Due to (11) we have
$$q_1 = g_t\mid_{t=0} = -2 u_{0_x}\mid_{t=0} = 4 \Lambda_{1_x} \neq
0,$$ so $w_{21} = 0$ and it means that $W_1 = 0.$

By induction, suppose that  $W_{k}(x,y) = 0$ for any $k < n.$
Substituting the expansion of $V$ into (14) we obtain the following
identities:
$$
n t^{n} q_1 w_{2n} + O(t^{n+1}) = 0,
$$
$$
n t^{n-1} w_{3n} + O(t^{n}) =0,
$$
$$
n t^{n-1} w_{4n} + O(t^{n}) =0,
$$
and by the same reasons we obtain $w_{2n} = w_{3n} = w_{4n} = 0.$
So we have $W_n=0.$ Lemma 2 is proved.

Consider now the metric which we have constructed. We shall use the
expansion (16) of $U$ to calculate several first terms of $\Lambda.$
Let's substitute the expansions into the system (11).  Calculating
the corresponding terms, one obtains the following identities:
$$
\lambda_1(x,y) = 0,
$$
$$
\lambda_2(x,y) = 2 ((\Lambda_1^{'}(x))^2 + (\Lambda_2^{'}(y))^2) +(\Lambda_1(x)+\Lambda_2(y)) (\Lambda_1^{''}(x) + \Lambda_2^{''}(y) ),
$$
it then follows that the metric $\Lambda$ is not Liouville class.

Notice that at the initial moment $t=0$ we have
$$
\frac{\partial \Omega}{\partial t} = \frac{1}{4} \frac{\partial}{\partial t} (g_x-f_y) =
-\frac{1}{2} (u_{0_{xx}}+u_{0_{yy}}) = \Lambda_1^{''}(x)-\Lambda_2^{''}(y) \neq 0.
$$
This means that at $t \neq 0$ the nonzero magnetic field appears.
Theorem 1 is proved.
\section{Quadratic integrals of magnetic flows}
In this section we prove the following:
\begin{theorem}
Consider the magnetic flow of the Riemannian metric
$ds^2=\Lambda(dx^2+dy^2)$ with the non-zero magnetic form $\Omega.$
Suppose the magnetic flow admits a first integral $F_2$ on all
energy levels such that $F_2$ is quadratic in momenta. Then in some
coordinates the functions $\Lambda$, $\Omega$ have the form of the
Example 1, so there exists another integral $F_1$ linear in momenta
and $F_2$ can be written as a combination of $H$ and $F_1.$
\end{theorem}
\medskip
{\it Proof.} Write the equation $\{F,H\}_{mg}=0$ for all points of
the phase space (not restricting to one energy level), where the
Hamiltonian $H$ and the integral $F_2$ are written as:
$$H=\frac{1}{2\Lambda}(p_1^2+p_2^2);\quad
F_2=a_{0}p_1^2+2a_{1}p_1p_2+a_2p_2^2+b_0p_1+b_1p_2+c_0.$$

Equating to zero coefficients of monomials $p_1^kp_2^l$ we get the
following 3 groups (a), (b), (c) of equations depending on the value
$k+l=3,2,1$ respectively. The first group (a) is:

$$\frac{(a_0)_x}{\Lambda}-a_0\left( \frac{1}{\Lambda} \right)_x-
a_1\left( \frac{1}{\Lambda} \right)_y=0, \eqno(a1)
$$
$$\frac{(a_0)_y}{\Lambda}+2\frac{(a_1)_x}{\Lambda}-
a_1\left( \frac{1}{\Lambda} \right)_x- a_2\left( \frac{1}{\Lambda}
\right)_y=0, \eqno(a2)
$$
$$\frac{(a_2)_x}{\Lambda}+2\frac{(a_1)_y}{\Lambda}-
a_0\left( \frac{1}{\Lambda} \right)_x- a_1\left( \frac{1}{\Lambda}
\right)_y=0, \eqno(a3)
$$
$$\frac{(a_2)_y}{\Lambda}-a_1\left( \frac{1}{\Lambda} \right)_x-
a_2\left( \frac{1}{\Lambda} \right)_y=0. \eqno(a4)
$$
The first group is the same as for geodesic flow without magnetic
field. Therefore in a usual way (following~\cite{12}) taking
(a1) -- (a3) and (a2) -- (a4) we arrive to the restrictions $a_1=const$
and $a_0-a_2=const.$ Without loss of generality one can assume
$a_1\equiv 0, a_0-a_2\equiv 4.$ Then substituting back to the (a1)
and (a4) we get
$$a_0=\frac{\phi(y)}{\Lambda},\quad
a_2=\frac{\psi(x)}{\Lambda},\quad 4\Lambda=\phi(y)-\psi(x).$$ The
next group of equations (b) reads:
$$\frac{(b_0)_x}{\Lambda}-b_0\left( \frac{1}{2\Lambda} \right)_x-
b_1\left( \frac{1}{2\Lambda} \right)_y=0, \eqno(b1)
$$
$$\frac{(b_1)_y}{\Lambda}-b_0\left( \frac{1}{2\Lambda} \right)_x-
b_1\left( \frac{1}{2\Lambda} \right)_y=0, \eqno(b2)$$
$$
\frac{(b_1)_x}{\Lambda}+\frac{(b_0)_y}{\Lambda}-\Omega\frac{2(a_2-a_0)}{\Lambda}=0.
\eqno(b3)
$$
Let us rename $b_0=2f,\ b_1=-2g.$ Then the equations (b1), (b2), (b3) can
be written in the form:
$$
f_x+\frac{\Lambda_x}{2\Lambda}f-\frac{\Lambda_y}{2\Lambda}g=0,
\eqno(b1')
$$
$$
f_x+g_y=0,\eqno(b2')
$$
$$
\Omega=\frac{1}{4}(g_x-f_y). \eqno(b3')
$$
Using $(b3')$ the equations of the last group (c) are :

$$c_{0_x}=-\frac{1}{2} g(g_x-f_y),\eqno(c1)$$

$$c_{0_y}=-\frac{1}{2} f(g_x-f_y). \eqno(c2)$$
One can easily see that using $(b2')$ one can rewrite (c1), (c2):
$$(c_0+\frac{g^2-f^2}{4})_x-\frac{1}{2}(fg)_y=0,\ \ (c_0+\frac{g^2-f^2}{4})_y+
\frac{1}{2}(fg)_x=0.$$ The last equations imply
$$
 c_0+\frac{g^2-f^2}{4}=K_1,\quad fg=K_2,
$$
where $K_1, K_2$ are constants.

Now multiplying $(b1')$ by $f\Lambda$ we get the following
$$
 (f^2\Lambda)_x=K_2\Lambda_y=K_2\frac{\phi'(y)}{4},
$$
and analogously for $g$:
$$(g^2\Lambda)_y=K_2\Lambda_x=-K_2\frac{\psi'(x)}{4}.$$
If the constant $K_2$ is not zero, then the last two equations imply
that the functions $\phi, \psi$ must be constant since otherwise
$f^2\Lambda$ or $g^2\Lambda$ would be unbounded. Thus $K_2=0.$ But
then we have that $$f^2\Lambda=F(y), \quad g^2\Lambda=G(x).$$
Multiplying these two and using $fg\equiv 0$ we conclude that
$$F(y)G(x)\equiv 0,$$
meaning that at least one say $F=f=0.$ Then from (b2') we have that
$$g=g(x).$$ And we have
$$g^2\Lambda=g^2(x)\frac{\phi(y)-\psi(x)}{4}=G(x).$$
Therefore the function $\phi(y)$ must be constant. Thus the
conformal factor $\Lambda$ must be a function of one variable
$\Lambda=\Lambda(x),$ and moreover the magnetic field is
$\Omega=\frac{1}{4} g'(x).$ In the Example 1 of the Introduction we
saw that this is exactly the case when the linear in momenta
integral $F_1$ exists and it is not hard to verify that the
quadratic integral $F_2$ can be expressed through $F_1$. This completes the proof.
\newpage

\renewcommand{\refname}{\normalfont\selectfont\large {\bf References}}

\end{document}